\documentclass{article}[12pt]
\usepackage{amsmath,amsfonts,amsthm}
\usepackage{amssymb}

\usepackage{epsfig}

\usepackage{subfigure}
\usepackage{graphics}
\usepackage{graphicx}
\usepackage{color}
\usepackage{soul}
\usepackage[makeroom]{cancel}

\usepackage{latexsym}
\usepackage{graphics}
\usepackage{amsmath}
\usepackage{amsthm}
\usepackage{xspace}
\usepackage{amssymb}
\usepackage{color}
\usepackage{cite}
\usepackage{latexsym}
\usepackage{graphics}
\usepackage{amsmath}
\usepackage{amsthm}
\usepackage{xspace}
\usepackage{amssymb}
\usepackage{epsfig}
\usepackage{comment}

\newcommand{\beql}[1]{\begin{equation}\label{#1}}
\newcommand{\eeql}{\end{equation}}
\newcommand{\eqn}[1]{(\ref{#1})}

\newcommand{\R}{\mathbb{R}}

\newcommand{\ci}{{\cal I}}
\newcommand{\cs}{{\cal S}}
\newcommand{\ck}{{\cal K}}
\newcommand{\cx}{{\cal X}}
\newcommand{\cm}{{\cal M}}

\newcommand{\ch}{{\cal H}}

\newcommand{\bk}{\boldsymbol{k}}
\newcommand{\bx}{\boldsymbol{x}}
\newcommand{\bX}{\boldsymbol{X}}
\newcommand{\bbeta}{\boldsymbol{\eta}}
\newcommand{\bnu}{\boldsymbol{\nu}}
\newcommand{\be}{\boldsymbol{e}}
\newcommand{\bu}{\boldsymbol{u}}
\newcommand{\bZero}{\boldsymbol{0}}

\newcommand{\ba}{\boldsymbol{a}}

\newtheorem{thm}{Theorem}
\newtheorem{lem}[thm]{Lemma}
\newtheorem{prop}[thm]{Proposition}

\newtheorem{assumption}[thm]{Assumption}
\newtheorem{definition}[thm]{Definition}
\newtheorem{conjecture}[thm]{Conjecture}

\theoremstyle{remark}
\newtheorem{remark}[thm]{Remark}

\begin{document}

\title{Large-scale heterogeneous
service systems with general packing constraints
}

\author
{
Alexander L. Stolyar \\
Lehigh University\\
200 W. Packer Ave., Room 484\\
Bethlehem, PA 18015 \\
\texttt{stolyar@lehigh.edu}
}

\date{\today}

\maketitle

\begin{abstract}

A service system with multiple  types of customers, arriving according to Poisson processes, is considered. 
The system is heterogeneous in that the servers also can be of
multiple types. 
Each customer has an independent exponentially distributed service time, 
with the mean determined by its type. 
Multiple customers (possibly of different types) can be placed for service into 
one server, subject to ``packing'' constraints, which depend on the server type.
Service times of different customers are independent, even if served simultaneously
by the same server.   
The large-scale asymptotic regime is considered
such that the customer arrival rates grow to infinity.

We consider two variants of the model. 
For the {\em infinite-server} model, we prove 
asymptotic optimality of the {\em Greedy Random} (GRAND) algorithm
in the sense
of minimizing the weighted (by type) number of occupied servers in steady-state.
(This version of GRAND generalizes that introduced in \cite{StZh2013} for the
homogeneous systems, with all servers of same type.)
We then introduce a natural extension of GRAND algorithm for 
{\em finite-server} systems with blocking.
Assuming subcritical system load, we 
prove existence, uniqueness, and local stability of the large-scale system equilibrium point
such that no blocking occurs.
This result strongly suggests a conjecture that the steady-state blocking probability 
under the algorithm vanishes in the large-scale limit.

\end{abstract}




\section{Introduction}
\label{sec-intro}

We consider a heterogeneous service system where servers can be of multiple types.
There are also multiple types of customers, each arriving according to an independent
Poisson process. 
Each customer has an independent exponentially distributed service time, 
with the mean determined by its type. 
Multiple customers (possibly of different types) can be placed for service into 
one server, subject to ``packing'' constraints, which depend on the server type.
Service times of different customers are independent, even if served simultaneously
by the same server.  Such a system arises, for example, as a model of dynamic real-time assignment
of virtual machines (``customers'') to physical host machines (``servers'') in a network cloud \cite{Gulati2012}, where typical objectives 
may be to minimize the number of occupied (non-idle) hosts or to minimize blocking/waiting of 
virtual machines. In this paper we consider two variants of the system, and study their properties
in the large-scale asymptotic regime, when the customer arrival rates (and then the number of occupied servers) are large.

The first variant of the system is such that there is an infinite ``supply'' of servers of each type. 
Each arriving customer
is assigned to a server immediately upon arrival. 
The asymptotic regime is considered such that
the customer arrival rates grow in proportion to a scaling parameter $r\to\infty$. 
Each server type $s$ is assigned a weight (``cost'') $\gamma_s$, and
the objective is to minimize the weighted number (``total cost'') of occupied servers in steady-state.
We prove that a generalized version of the {\em Greedy Random} (GRAND) algorithm, introduced
in \cite{StZh2013} for a homogeneous system (with one server type), is asymptotically optimal,
in the sense described below in this paragraph.
The basic idea of GRAND is to assign an arriving customer of a given type $i$ to a server
chosen randomly uniformly among servers available to it, i.e. those servers where a type $i$ customer
can be added without violating packing constraints. A particular GRAND algorithm that we consider
for the infinite server system, which is labeled GRAND($\ba Z$), 
is as follows. There is a parameter $a_s>0$ for each server type $s$; $\ba = (a_s)$ is the vector with components $a_s$. An arriving customer
picks uniformly at random an available server among all currently occupied servers plus designated 
numbers $a_s Z$ of idle servers (called ``zero-servers'') of each type $s$, where $Z$ is the current total number of all customers. (GRAND($a Z$) algorithm of \cite{StZh2013} is a special case of GRAND($\ba Z$),
with single parameter $a>0$, because there is only one server type.) {\em GRAND($\ba Z$) achieves optimality
if we first take the limit of system stationary distributions as $r\to\infty$, 
and then take the limit on $a_s= \alpha^{\gamma_s}\downarrow 0$,
with common parameter $\alpha\downarrow 0$.} (We believe that a stronger form of asymptotic 
optimality, when only the limit $r\to\infty$ is taken, holds for a different version of  
GRAND, with the number
of zero-servers of type $s$ equal to $Z^{(p-1)\gamma_s+1}$, where parameter $p<1$ is close to $1$.
See Conjecture~\ref{conj-zp} at the end of
Section~\ref{sec-results-infinite}.)

It is important to emphasize that GRAND($\ba Z$) achieves asymptotic optimality {\em without 
utilizing any knowledge of the system structural parameters}. Namely, the algorithm need not ``know" the server
types or exact states of the currently occupied servers. All it needs to know about each currently occupied server
is whether or not it can ``accept'' an additional customer of type $i$, for each $i$.
Note that the setting of the algorithm parameters $a_s$, that achieves asymptotic optimality, 
depends only on the weights $\gamma_s$,
which are the parameters of the objective (as opposed to system parameters).
One of the key qualitative insights of \cite{StZh2013} was the surprising fact that an algorithm 
as simple as GRAND can be asymptotically optimal. The fact that an appropriately
generalized, but still extremely simple, version of GRAND is optimal for in a heterogeneous system,
is still more surprising.

The second variant is a system with finite size pools of servers of each type. 
Each arriving customer can be either immediately assigned to a server or immediately blocked (in which case
it leaves the system without receiving service).
The asymptotic regime is such that both the arrival rates and the server pool sizes scale in proportion to parameter $r\to\infty$.
We consider a different version of the GRAND algorithm, labeled GRAND-F, which simply assigns 
each arriving customer randomly uniformly to any available to it server in the system,
and blocks the customer if there are no such available servers. 
We study the dynamics of the fluid paths (obtained by ``fluid" scaling and then the 
$r\to\infty$ limit). Assuming the system is subcritically loaded, we prove existence,
uniqueness and local stability of a system equilibrium point, such that there is no blocking.
These results strongly suggest
a conjecture that GRAND-F is asymptotically optimal in that,
under subcritical load, the limit
of the system stationary distributions is concentrated 
on the equilibrium point described above, and therefore {\em the steady-state blocking probability vanishes
in the $r\to\infty$ limit}.
We note that the equilibrium point local stability property
is stronger than a typical
``fixed point'' argument, based on the assumption of asymptotic independence of server states
 (or, ``independence ansatz,'' in the terminology of \cite{BLP2012-jsq-asymp-indep,BLP2013-jsq-asymp-tail}). 
The fixed point argument allows one to characterize (and then possibly derive)
the limit of the stationary distributions, assuming the ansatz holds. If the ansatz is proved, this of course proves
the limit of the stationary distributions. If the ansatz is {\em not}
proved, the fixed point argument is equivalent to the property that the equilibrium point 
is an invariant point of the fluid paths. 
The local stability of the equilibrium point that we prove,
is a stronger property than just its existence and invariance,
and therefore it provides a stronger support for the asymptotic optimality conjecture.
(The relation between the local stability and the fixed point argument is discussed in detail in Section~\ref{sec-fixed-point}.)

We want to emphasize that the packing constraints that we consider are extremely
general. (They are of the same kind as those in \cite{St2012,StZh2012,StZh2013}; we additionally allow
them to depend on the server type.) In particular, they are far more general than {\em vector packing} constraints.
Vector packing refers to the situation when a server has the corresponding resource-vector,
giving the amounts of resources of different types that it possesses;
for each customer type there is the requirement-vector, giving the resource requirements of one customer;
the constraint is that the sum of the requirement-vectors of the customers placed into a server cannot exceed
its resource-vector. Packing of virtual machines into physical machines in a network cloud \cite{Gulati2012} is 
an example of vector packing. 

Finally, we note that GRAND-F can be very efficiently implemented via a ``pull-based'' mechanism 
(see \cite{St2014_pull} and references therein), which has a very low signaling message exchange rate 
between the ``router'' and the servers. In fact,
GRAND-F algorithm can be viewed as an extension of PULL algorithm \cite{St2014_pull}
to service systems with packing constraints. 
(This is discussed in more detail in Remark~\ref{rem-pull} in Section~\ref{sec-results-finite}.)

\subsection{Related previous work}
\label{sec-rel-prev-work}

As mentioned above, the main practical motivation for our model is the problem of real-time dynamic assignment
of virtual machines (VM) to physical host machines (PM) in a network cloud. (A general discussion of the issues that arise
in this application can be found in \cite{Gulati2012}.) Since multiple VMs can simultaneously occupy 
(be ``packed into'') same PM, this naturally leads to bin packing type models.
There is an extensive literature on the classical bin packing (see, e.g., \cite{Csirik2006,Bansal2009,GR2012}
for reviews and recent results), where each ``item'' (customer) once placed into a ``bin'' (server) stays in that bin forever. However, the dynamic VM-to-PM assignment problem is such that each VM (customer) 
leaves its PM (server), and the system, after its service is completed. This in turn naturally leads 
the models that we consider, i.e. service systems with packing constraints at the servers.

The infinite-server variant of our model is a generalization of the homogeneous (one server type)
model studied in \cite{St2012,StZh2012,StZh2013}, which focused on the problem of minimizing the number 
of occupied servers in steady-state. In particular, GRAND algorithm was proposed and shown to be asymptotically optimal in \cite{StZh2013}. (Papers \cite{St2012,StZh2012} have studied a different algorithm,
which needs to know the structure of packing constraints and to use the exact current states of all servers.)
Our model allows, in addition, multiple server types
and we consider a more general problem of minimizing the weighted number of servers; the analysis
of this variant of our  model is a generalization of that in \cite{StZh2013}. 
A homogeneous infinite-server model, specialized to vector packing constraints,
was also considered in \cite{GZS2013}, where a randomized version of Best Fit algorithm was proved asymptotically optimal.

The finite-server variant of our model is related to 
the model in recent paper 
\cite{XDLS}, which considers blocking in a
homogeneous system, specialized to one-dimensional (single resource) vector packing constraints.
(In \cite{XDLS} all servers
are of the same type, and
the term {\em heterogeneous} refers to multiple customer types, 
which our model also allows. So, in our terminology, the system in \cite{XDLS} is homogeneous.)
The algorithm in \cite{XDLS} is of the {\em power-of-d-choices} type 
\cite{VDK96,Mitz2001,BLP2012-jsq-asymp-indep,BLP2013-jsq-asymp-tail}, 
namely each arriving
customer goes to the server which has the largest amount of unused resource, out of the $d$
servers chosen uniformly at random. 
The paper
uses a fixed point argument (independence ansatz) to derive the form of the equilibrium point,
which is conjectured to be the asymptotic limit of the system steady-state.
(In addition, the paper derives some performance bounds.) 
Of course, the equilibrium point  under the power-of-d-choices algorithm is different from that under our GRAND-F algorithm.
It is such that the blocking probability does {\em not} (and cannot be expected to)
vanish in the limit. Therefore, the relation between the power-of-d-choices algorithm and GRAND-F for the systems
with packing constraints, is analogous to the relation between power-of-d-choices and PULL algorithm 
\cite{St2014_pull} for service systems
without packing, where the blocking (or waiting) probability vanishes under PULL, but not under the power-of-d-choices. (GRAND-F can be viewed as an extension of PULL algorithm to systems with packing constraints.
See Remark~\ref{rem-pull} in Section~\ref{sec-results-finite}.)

Papers \cite{Maguluri2012, Maguluri2013} consider a homogeneous finite-server system with queues (and no blocking), and focus on the system stability (or, throughput maximization).
In \cite{GSW2012} a heterogeneous finite-server system is considered, with the objective of minimizing 
maximum load across server pools; the algorithms proposed in \cite{GSW2012} essentially treat the system as an infinite-server one. The algorithms in \cite{Maguluri2012, Maguluri2013,GSW2012} 
are completely different from the variants of GRAND algorithm studied in this paper.

\subsection{Layout of the rest of the paper}

Basic notation used throughout the paper is given in Section~\ref{subsec-notation}.
The model and the main results are stated in Section~\ref{sec-model}.
The basic structure of the system, common to both variants, is given in Section~\ref{sec-packing-constr}. The infinite-server system, GRAND($\ba Z$) algorithm and the main results for it (Theorems~\ref{th-grand-fluid} and \ref{th-grand-fluid-convergence}) are presented in Section~\ref{sec-results-infinite}.
Section~\ref{sec-results-finite} defines the finite-server system, GRAND-F algorithm, and states the main result for it
informally in Proposition~\ref{prop-main-finite} (with formal statements given later in 
Lemmas~\ref{lem-fsp-conv-local-on-diamond} and \ref{lem-fsp-conv-local}). Sections~\ref{sec-opt-conv}
and \ref{sec-fsp-dynamics} 
 contain proofs of the infinite-server/GRAND($\ba Z$) results,
while Section~\ref{sec-fsp-dynamics-finite} contain those for finite-server/GRAND-F. Concluding remarks
are given in Section~\ref{sec-further-work}.

\subsection{Basic notation}
\label{subsec-notation}

Sets of real and real non-negative numbers are denoted by $\R$ and $\R_+$, respectively.
 We use bold and plain letters  for vectors and scalars, respectively.
The standard Euclidean norm of a vector $\bx\in \R^n$ is denoted by $\|\bx\|$. 
Convergence $\bx \to \bu \in \R^n$ means ordinary convergence in $\R^n$,
while $\bx \to U \subseteq \R^n$ means convergence to a set, namely,
$\inf_{\bu\in U} \|\bx-\bu\|\to 0$.
The $i$-th coordinate unit vector in $\R^n$ is denoted by $\be_i$.
Symbol $\implies$
denotes convergence in distribution of random variables taking values in space $\R^n$
equipped with the Borel $\sigma$-algebra. The
abbreviation {\em w.p.1} means convergence {\em with probability 1}.
We often write $x(\cdot)$ to mean the function (or random process) $\{x(t),~t\ge 0\}$.
Abbreviation {\em u.o.c.} means 
{\em uniform on compact sets} convergence of functions.
The cardinality of a finite set $\mathcal{N}$ is 
$|\mathcal{N}|$. Indicator function $I\{A\}$ for a condition $A$ is equal to $1$
if $A$ holds and $0$ otherwise.  $\lceil \xi \rceil$ denotes the smallest integer
greater than or equal to $\xi$, and $\lfloor \xi \rfloor$ denotes the largest integer 
smaller than or equal to $\xi$.
For a finite set of scalar functions $f_n(t), ~t\ge 0$, $n\in\mathcal{N}$, a point $t$ is called
{\em regular} if for any subset $\mathcal{N}' \subseteq \mathcal{N}$ the 
derivatives
$\frac{d}{dt} \max_{n\in\mathcal{N}'} f_n(t)$ and
$\frac{d}{dt} \min_{n\in\mathcal{N}'} f_n(t)$
exist.

\section{Model and main results}
\label{sec-model}

In this section we formally define the two variants of the model with heterogeneous servers, and state our main results for them. 
The first variant is a generalization of the infinite-server model in \cite{St2012,StZh2012,StZh2013} in that 
we allow different types of servers, as opposed to just one type.
The number of servers of each type is infinite and there is no blocking of arriving customers. For this version of the model the underlying objective is to minimize the weighted number of occupied servers in steady-state.
The second variant is the model with different server types, but with
finite number of servers of each type. If an arriving customer cannot be immediately assigned to some server in the system, it is blocked. In such a system, the underlying objective is to minimize blocking. Before defining these two variants of the model, in the next subsection
we define the basic structure of the system
(most importantly the server packing constraints),
which is common for both model variants.

\subsection{Heterogeneous servers. Packing constraints}
\label{sec-packing-constr}

We consider a service system with $I$  types 
of customers, indexed by $i \in \{1,2,\ldots,I\} \equiv \ci$. 
The service time of a type-$i$ customer is an exponentially distributed random variable with mean $1/\mu_i$.
All  customers' service times are mutually independent.
There are $S$ types of servers, indexed  $s \in \{1,2,\ldots,S\} \equiv \cs$, and infinite ``supply'' of servers of each type.
A server of each type can potentially serve more than one customer simultaneously, subject to the following very general packing constraints. We say that a vector $\bk = (k_1,\ldots,k_I;s)$ with non-negative integer
$k_i, ~i\in \ci,$ and $s\in \cs$ is a server {\em configuration}, if a type $s$ server can simultaneously serve a combination of customers of different types given by the values $k_i$. 
A configuration $\bk$ with specific value of $s$ is a type $s$ server configuration.
For any $s$, there is a finite set of all allowed type $s$ server configurations, denoted by $\bar\ck^s$.
We assume that $\bar\ck^s$ satisfies a natural {\em monotonicity} condition: if 
$\bk\in \bar\ck^s$, then all ``smaller'' configurations $\bk'= (k'_1,\ldots,k'_I;s)$, i.e. such that $k'_i \le k_i$ for all $i$, belong to $\bar\ck^s$ as well. Without loss of generality, assume that for each $i$, $(\be_i;s) \in \bar\ck^s$ for at least one $s$, where $\be_i$ is the $i$-th coordinate unit vector (otherwise, type-$i$ customers cannot be served at all).
By convention, for any $s$, vector $\bZero^s \equiv (\bZero;s) \in \bar\ck^s$, where
$\bk=\bZero$ is the $I$-dimensional component-wise zero vector -- this is the configuration of an empty 
type $s$ server. We denote by $\ck^s=\bar\ck^s \setminus \{\bZero^s\}$ 
the set of type $s$ server configurations {\em not} including the empty (or, zero) configuration. 
Denote by $\bar\ck = \cup_s \bar\ck^s$ and $\ck = \cup_s \ck^s$ the sets of all configurations and all non-zero
configurations, respectively. In what follows, we use the following slightly abusive notations: for $\bk\in \bar\ck$,
$\bk+\be_i$ means vector $\bk$ with $k_i$ replaced by $k_i+1$, and similarly for $\bk-\be_i$.

An important feature of the model is that simultaneous service does {\em not} affect the service time
distributions of individual customers. In other words, the service time of a customer is unaffected by whether or not there are other customers served simultaneously by the same server. A customer can be ``added'' to an empty or occupied server, as long as the packing constraints are not violated.
Namely, a type $i$ customer can be added to a server of type $s$ whose current configuration $\bk\in\bar\ck^s$  is such that $\bk+\be_i \in \ck^s$. 
When the service of a type-$i$ customer by a  server in configuration $\bk$ is completed,
the customer leaves the system and the server's configuration changes to $\bk-\be_i$.

\subsection{Infinite-server system}
\label{sec-results-infinite}

In this section we define the infinite-server system, the proposed generalized GRAND($aZ$) assignment (or packing) algorithm, and state the asymptotic optimality results for this algorithm.

We consider a system, as described in Section~\ref{sec-packing-constr},
in which there is an infinite ``supply'' of servers of each type $s\in\cs$.
Customers of type $i$ arrive as an independent Poisson process of rate $\Lambda_i >0$;
these arrival processes are independent of each other and of the customer service times.
Each arriving customer is immediately placed for service in one of the servers, as long as packing 
constraints are not violated.

Denote by $X_{\bk}$ the number of servers in configuration $\bk\in \ck^s$. The system state is then the vector $\bX = \{X_{\bk}, ~\bk\in \ck\}$. 

A {\em placement algorithm} (or packing rule) determines where an arriving customer is placed, as a function of the current system state $\bX$. Under any well-defined placement algorithm, 
the process $\{\bX(t), t\ge 0\}$ is a continuous-time Markov chain 
with a countable state space. It is easily seen to be irreducible and positive recurrent: the positive recurrence follows from the fact
that the total number $Y_i(t)$ of type-$i$ customers in the system is independent from the  placement algorithm, and its stationary distribution is 
Poisson with mean $\Lambda_i/\mu_i$; we denote by $Y_i(\infty)$ the random value of $Y_i(t)$ in steady-state -- it is, therefore, a Poisson random 
variable with mean $\Lambda_i/\mu_i$.
Consequently, the process $\{\bX(t), ~t\ge 0\}$ has a unique stationary distribution; let $\bX(\infty) = \{X_{\bk}(\infty), \bk \in \ck \}$ be the random system state $\bX(t)$ in stationary regime.

We are interested in finding a 
placement algorithm that minimizes 
the total weighted number of  occupied servers $\sum_{\bk\in\ck} X_{\bk}(\infty)$
in the stationary regime. 

Consider the following generalization of the Greedy-Random (GRAND) algorithm, introduced in \cite{StZh2013}. More specifically, it is a generalization of 
the special form of the algorithm, called in \cite{StZh2013} GRAND($aZ$).

\begin{definition}[Greedy-Random (GRAND($\ba Z$)) algorithm for heterogeneous infinite-server systems]\label{df:grand}
The algorithm is parameterized by a vector $\ba=(a_s, ~s\in \cs)$ of real numbers $a_s>0$.
Let $Z(t)=\sum_i \sum_{\bk} k_i X_{\bk}(t)$ denote the total number of customers in the system at time $t$.
At any given time $t$, there is a designated finite set of $X_{\bZero^s}(t) =\lceil a_s Z(t) \rceil\ge 0$ empty type $s$ servers, 
called {\em $s$-zero-servers}.  \\
A new customer, 
say of type $i$, arriving at time $t$
is placed into a server chosen randomly
 uniformly among those zero-servers (of any type $s$) and occupied servers,
where it can still fit.
 In other words,
the total number of servers available to a type-$i$ arrival at time $t$ is
$$
X_{(i)}(t) \doteq  \sum_{\bk\in \bar \ck:~\bk+\be_i\in \ck} X_{\bk}(t) \equiv 
\sum_{s:~\be_i \in \ck^s} \left[X_{\bZero^s}(t) + \sum_{\bk\in \ck:~\bk+\be_i\in \ck} X_{\bk}(t)\right].
$$
If $X_{(i)}(t)=0$, the customer is placed into an empty server of any type $s$ such that $\be_i \in \ck^s$.
\end{definition}

The GRAND($\ba Z$) algorithm is easily implementable. 
(A detailed discussion of the implementation issues of the GRAND algorithm is given below in
Remark~\ref{rem-pull}, in the context of finite-server systems.)

We now define the asymptotic regime. Let $r\to\infty$ be a positive scaling parameter.
 More specifically, assume that $r\ge 1$, and $r$ increases to infinity along a discrete
sequence. 
Customer arrival rates scale linearly with $r$; namely,
for each $r$, $\Lambda_i = \lambda_i r$, where $\lambda_i$ are fixed positive parameters.
Let $(\bX^r(t), ~t\ge 0),$ be the process associated with a system with parameter $r$, and 
let $\bX^r(\infty)$ be the (random) system state in the stationary regime.
(Note that we do {\em not} include the zero-server numbers $X_{\bZero^s}^r(t)$ into $\bX^r(t)=\{X_{\bk}^r(t),~\bk\in \ck\}$.)
For each $i$, denote by $Y^r_i(t) \equiv \sum_{\bk\in\ck} k_i X^r_{\bk}(t)$ the total number
of customers of type $i$. Since arriving customers are placed for service immediately
and their service times are independent of each other and of the rest of the system, 
 $Y^r_i(\infty)$ is a Poisson random variable with mean $r \rho_i$, where $\rho_i\equiv \lambda_i/\mu_i$.
Moreover, $Y^r_i(\infty)$ are independent across $i$.
Since the total number of occupied servers is no greater than the total number 
of 
customers, 
$\sum_{\bk\in\ck} X_{\bk}^r(t) \le Z^r(t)\equiv \sum_i Y^r_i(t)$, we have a simple upper bound 
on the total number of occupied servers in steady state,
$\sum_{\bk\in\ck}  X_{\bk}^r(\infty) \le Z^r(\infty) \equiv \sum_i Y^r_i(\infty)$,
where $Z^r(\infty)$ is a Poisson random variable with mean $r \sum_i \rho_i$.
Without loss of generality, from now on  we assume $\sum_i \rho_i=1$.
This is equivalent to rechoosing the parameter $r$ to be $r \sum_i \rho_i$.

The {\em fluid-scaled} process is $\bx^r(t)=\bX^r(t)/r$, $t \in [0, \infty)$. 
We also define $\bx^r(\infty) = \bX^r(\infty)/r$.
For any $r$, $\bx^r(t)$ takes values in the  non-negative orthant $\R_+^{|\ck|}$.
Similarly, $y^r_i(t)=Y^r_i(t)/r$, $z^r(t)=Z^r(t)/r$, $x^r_{\bZero^s}(t)=X^r_{\bZero^s}(t)/r$ and 
$x^r_{(i)}(t)=X^r_{(i)}(t)/r$,  for  $t \ge 0$ and $t=\infty$.
Since $\sum_{\bk\in \ck}  x_{\bk}^r(\infty) \le z^r(\infty)=Z^r(\infty)/r$,
we see that the random variables
$(\sum_{\bk\in \ck} x_{\bk}^r(\infty))$ are uniformly integrable in $r$.
This in particular implies that the sequence of distributions of $\bx^r(\infty)$ is tight,
and therefore there always exists a limit $\bx(\infty)$ in distribution, so that $\bx^r(\infty)\implies \bx(\infty)$,
along a subsequence of $r$.

The limit (random) vector $\bx(\infty)$ satisfies the following conservation laws:
\beql{eq-cons-laws}
\sum_{\bk\in\ck} k_i x_{\bk}(\infty) \equiv y_i(\infty) = \rho_i, ~~\forall i,
\end{equation}
and, in particular, 
\beql{eq-cons-laws2}
z_i(\infty)\equiv \sum_i y_i(\infty) \equiv \sum_i \rho_i
= 1.
\end{equation}
Therefore, the values of $\bx(\infty)$ are confined to the convex compact 
$(|\ck|-I)$-dimensional 
polyhedron
$$
\cx \equiv \{\bx\in \R_+^{|\ck|} ~|~  \sum_s \sum_{\bk\in \ck^s} k_i x_{\bk} = \rho_i, ~\forall i\in\ci \}.
$$
We will slightly abuse notation by using symbol $\bx$ for a generic element of $\cx$;
while $\bx(\infty)$, and later $\bx(t)$, refer to random  elements taking values in $\cx$.

Also note that under GRAND($aZ$), for any server type $s$, $x^r_{\bZero^s}(\infty) \implies x_{\bZero^s}(\infty)=a_s z(\infty)=a_s$, as $r\rightarrow \infty$.

The asymptotic regime and the associated basic properties \eqn{eq-cons-laws}
and \eqn{eq-cons-laws2} hold {\em for any placement algorithm}. 
 Indeed, \eqn{eq-cons-laws}
and \eqn{eq-cons-laws2} only depend on the already mentioned  fact 
that all $Y_i^r(\infty)$ are mutually independent
Poisson  random variables
with means $\rho_i r$.

Let the server weights $\gamma_s>0$, $s\in \cs$, be fixed.
(One can think of $\gamma_s$ as 
the ``cost" rate of using one type $s$ server.)
Consider the following problem of minimizing the weighted number of occupied servers,
on the fluid scale:
$\min_{\bx\in \cx} \sum_{s\in\cs} \sum_{\bk\in\ck^s} \gamma_s x_{\bk}$. It is a
linear program:
\beql{eq-opt}
\min_{\bx\in\R_+^{|\ck|}} \sum_{s\in\cs} \sum_{\bk\in\ck^s} \gamma_s x_{\bk},
\end{equation}
subject to
\beql{eq-cons-laws222}
\sum_{\bk\in\ck} k_i x_{\bk} = \rho_i, ~~\forall i.
\end{equation}
Without loss of generality, assume that the weights are scaled so that $\gamma_1=1$.
Denote by $\cx^* \subseteq \cx$ the set of optimal solutions of \eqn{eq-opt}-\eqn{eq-cons-laws222}.

For future reference, we record the following observations and notation. 
Using the monotonicity of $\bar{\ck}$, it is easy to check that if in the LP 
\eqn{eq-opt}-\eqn{eq-cons-laws222} we replace equality constraints
\eqn{eq-cons-laws222} with the inequality constraints
\beql{eq-cons-laws222ge}
\sum_{\bk\in\ck} k_i x_{\bk} \ge \rho_i, ~~\forall i,
\end{equation}
the new LP \eqn{eq-opt},\eqn{eq-cons-laws222ge} has same optimal value, and its set of the optimal solutions 
$\cx^{**}$ contains $\cx^*$, or more precisely, $\cx^* = \cx^{**} \cap \cx$.
From here,  using Kuhn-Tucker theorem,
$\bx \in \cx^*$ if and only if there exists a vector $\bbeta=\{\eta_i, ~i\in \ci\}$
of Lagrange multipliers, corresponding to the
inequality  constraints \eqn{eq-cons-laws222ge}, such that the following conditions hold:
\beql{eq-dual-1}
\bx\in \cx,
\end{equation}
\beql{eq-dual-111}
\eta_i \ge 0, ~~\forall i \in \ci,
\end{equation}
\beql{eq-dual-2}
\sum_i k_i \eta_i \le \gamma_s, ~~\bk \in \ck^s,
\end{equation}
\beql{eq-dual-3}
\mbox{for $\bk \in \ck^s$, ~~ condition} ~\sum_i k_i \eta_i < \gamma_s ~\mbox{implies}~x_{\bk}=0.
\end{equation}
Vectors $\bbeta$  satisfying \eqn{eq-dual-1}-\eqn{eq-dual-3} for some $\bx \in \cx$
 are optimal solutions to the problem dual to LP \eqn{eq-opt},\eqn{eq-cons-laws222ge}.
 They form a convex set, which we denote by $\ch^*$; it is easy to check that $\ch^*$ is compact.

For each parameter-vector $\ba$ (as in the definition of GRAND($\ba Z$) algorithm), denote
\beql{eq-L-def}
L^{(\ba)}(\bx) = 
\sum_s \sum_{\bk\in\ck^s} x_{\bk} \log [x_{\bk} c_{\bk} /(e a_s)],
\end{equation}
where $c_{\bk} \doteq \prod_i k_i !$, $0!=1$. 
Then for $\bk\in \ck^s$ we have 
\beql{eq-L-partial}
(\partial / \partial x_{\bk}) L^{(\ba)}(\bx) =
\log [x_{\bk} c_{\bk} / a_s ].
\end{equation}
Note that if we adopt a convention that
\beql{eq-formal-deriv-zero}
(\partial / \partial x_{\bZero^s}) L^{(\ba)}(\bx)|_{x_{\bZero^s}=a_s} = 0,
\end{equation}
then \eqn{eq-L-partial} is valid for $\bk=\bZero^s$ and $x_{\bZero^s}=a_s$, which will be useful later.

The function $L^{(\ba)}(\bx)$ is strictly convex in $\bx\in\R_+^{|\ck|}$.
Consider the problem $\min_{\bx\in \cx} L^{(\ba)}(\bx)$. It is the
following convex optimization problem:
\beql{eq-opt-grand}
\min_{\bx\in\R_+^{|\ck|}} L^{(\ba)}(\bx),
\end{equation}
subject to
\beql{eq-cons-laws-grand}
\sum_{\bk\in\ck} k_i x_{\bk} = \rho_i, ~~\forall i.
\end{equation}
Denote by $\bx^{*,\ba} \in \cx$ its unique optimal solution. 
Using \eqn{eq-L-partial} it is easy to check that $x^{*,\ba}_{\bk}>0$ for all ${\bk} \in \ck$.
There exists a vector $\bnu^{*,\ba} = \{\nu_i^{*,\ba}, ~ i\in \ci\}$ of  Lagrange multipliers for 
the constraints \eqn{eq-cons-laws-grand}, such that $\bx^{*,\ba}$ solves problem
$$
\min_{\bx\in \R_+^{|\ck|}} L^{(\ba)}(\bx) + \sum_i \nu_i^{*,\ba} (\rho_i - \sum_{\bk\in\ck} k_i x_{\bk}).
$$
We see that $\log [x_{\bk}^{*,\ba} c_{\bk} / a_s ] - \sum_i \nu_i^{*,\ba} k_i = 0$, $ \bk \in \ck$. 
Therefore, $\bx^{*,\ba}$ has the product form
\beql{eq-grand-product}
x^{*,\ba}_{\bk} = \frac{a_s}{c_{\bk}} \exp \left[\sum_i k_i \nu_i^{*,\ba}\right], ~~ \bk \in \ck^s.
\end{equation}

This in particular implies that the Lagrange multipliers $\nu_i^{*,\ba}$ are unique 
and are equal to\\
 $\nu_i^{*,\ba}=\log (x_{\be_i}^{*,\ba} / a_s)$, by considering \eqref{eq-grand-product} for $\be_i$, $i \in \ci$; note also that they
can have any sign (not necessarily non-negative). Therefore, we obtain the following fact. {\em
A point $\bx\in \cx$ is the optimal solution to \eqn{eq-opt-grand}-\eqn{eq-cons-laws-grand}
(that is $\bx=\bx^{*,\ba}$) if and only if it has a product form representation
\eqn{eq-grand-product} for some vector $\bnu^{*,\ba}$.} (The 'only if' part we just proved, and the 'if' follows
from Kuhn-Tucker theorem.)

Our main results on the asymptotic optimality of GRAND($aZ$) algorithm
for the system with infinite number of servers
 are the following Theorems~\ref{th-grand-fluid} and \ref{th-grand-fluid-convergence}. 

\begin{thm}
\label{th-grand-fluid}
Let the parameter vector $\ba$ be fixed.
Consider a sequence of systems under the GRAND($\ba Z$) algorithm, indexed by $r$,
and let $\bx^r(\infty)$ denote the random state of the fluid-scaled process
in the stationary regime.
Then, as $r\to\infty$,
$$
\bx^r(\infty) \implies \bx^{*,\ba}.
$$
\end{thm}

\begin{thm}
\label{th-grand-fluid-convergence}
Suppose the parameter vector $\ba$ itself depends on a single parameter $\alpha>0$ as follows:
$a_s = \alpha^{\gamma_s}, s\in \cs$. Then,
as $\alpha\downarrow 0$, $\bx^{*,\ba} \to \cx^*$ and $(-\log \alpha)^{-1} \bnu^{*,\ba} \to \ch^*$.
\end{thm}

Theorems~\ref{th-grand-fluid} and \ref{th-grand-fluid-convergence} show that GRAND($\ba Z$) 
is asymptotically optimal in the sense that $\bx^r(\infty)$ converges to the optimal set $\cx^*$,
if we first take the limit $r\to\infty$, and then take the limit $\alpha\downarrow 0$ with $a_s = \alpha^{\gamma_s}$.

It was proved in a recent paper \cite{StZh2015} (which is posterior to this paper) that 
a stronger form of asymptotic optimality, when only the limit $r\to\infty$ is taken,
is achieved by the following version of GRAND, called GRAND($Z^p$). This is a GRAND algorithm
with the number of zero-servers depending on $Z$ as $Z^p$, where $p<1$ is a parameter, which is 
sufficiently close to $1$, but depends only on the packing constraints. GRAND($Z^p$) can be
informally interpreted as GRAND($aZ$), with $a$ being variable $a=Z^{p-1}$ rather than constant.
This suggests that for the heterogeneous infinite-server system that we consider, 
the stronger form of asymptotic optimality should hold, if we make $a_s$ variable,
equal to $Z^{(p-1)\gamma_s}$. Specifically, we believe that the methods of \cite{StZh2015}
can be extended to prove the following fact.

\begin{conjecture}
\label{conj-zp}
Consider the GRAND algorithm with the number
of zero-servers of type $s$ equal to $Z^{(p-1)\gamma_s+1}$, where parameter $p<1$ is 
sufficiently close to $1$, but depends only on the packing constraints (i.e., sets $\ck^s$).
Then, as $r\to\infty$, $d(\bx^r(\infty),\cx^*) \Rightarrow 0$, where $d(\bx,U)$ is the
distance from point $\bx$ to set $U$.
\end{conjecture}

\subsection{Finite-server system}
\label{sec-results-finite}

We now consider a version of the system, where the number of servers of each type is finite.
Namely, there is 
a finite number $H_s>0$ of servers of type $s$.
Customers of type $i$ arrive as an independent Poisson process of rate $\Lambda_i >0$
(and these processes are independent from the customer service times). 
Each arriving type $i$ customer can be either immediately placed for service into one of the servers (subject to
packing constraints) or immediately blocked, in which case it  immediately leaves the system. If there is no server where an arriving customer can be placed,
the customer is necessarily blocked. 

Let $X_{\bk}$ denote the number of servers in configuration $\bk\in \ck^s$ and the system state is the vector $\bX = \{X_{\bk}, ~\bk\in \ck\}$. (Same notation as for the infinite-server system.)
Note that we do not include the numbers $X_{\bZero^s}$ of
empty servers of each type (i.e., $s$-zero-servers) into the state $\bX$. However, those number are, of course, 
uniquely determined by $\bX$, because at all times we have the conservation law
$$
X_{\bZero^s} + \sum_{\bk\in\ck^s} X_{\bk} = \sum_{\bk\in\bar\ck^s} X_{\bk} = H_s, ~s\in\cs.
$$

In such a system, a  placement algorithm (or packing rule) determines, depending on the current system state $\bX$, whether or not an arriving customer is accepted (i.e., not blocked), and if so, into which server it is placed.
(If there are no servers, where a customer can be placed, it is necessarily blocked.)
Under any well-defined placement algorithm, 
the process $\{\bX(t), t\ge 0\}$ is a continuous-time Markov chain 
with finite state space; it is easily seen to be irreducible and, therefore, ergodic, with 
unique stationary distribution. Let $\bX(\infty) = \{X_{\bk}(\infty), \bk \in \ck \}$ be the random system state $\bX(t)$ in stationary regime. It is also easy to see that $Y_i(\infty)$ -- the steady-state random number
of all type $i$ customers in the system -- is stochastically dominated by 
that in the infinite-server system, i.e. by
a Poisson random 
variable with mean $\Lambda_i/\mu_i$.

For this system, the underlying objective is to minimize blocking in steady-state.
We consider the following version of the Greedy-Random (GRAND) algorithm, 
for the finite-server systems. It will be labeled GRAND-F. 

\begin{definition}[GRAND-F]\label{df:grand-b}
A new customer, 
say of type $i$, arriving at time $t$
is placed into a server chosen randomly
 uniformly among all servers in the system
where it can still fit.
(The total number of servers available for a type $i$ customer addition at time $t$ is
$$
X_{(i)}(t) \doteq  \sum_{\bk\in \bar \ck:~\bk+\be_i\in \ck} X_{\bk}(t). ~)
$$
If there are no such available servers (i.e., $X_{(i)}(t)=0$), the customer is blocked.
\end{definition}

\begin{remark}
\label{rem-pull}
An implementation of GRAND-F algorithm only requires that the ``router'' (an entity, making an assignment
decision for each arriving customer) knows which servers are currently available for an addition 
of a type $i$ customer, for each $i\in\ci$. The router does {\em not} need to know the exact configurations
of the servers. Moreover, it does {\em not} even need to know the server types!
Therefore, the router needs to maintain only $I$ bits of information
for each server. This in turn is easily achievable, for example, by using a {\em pull-based}
mechanism, analogous to that used by the PULL algorithm proposed in \cite{St2014_pull} (in a different context, for systems without non-trivial packing constraints).
A specific pull-based mechanism 
to work in conjunction with GRAND-F can be as follows.\\
(a) Upon a customer, say of type $i$, arrival,
the router follows GRAND-F rule for choosing a server. If there are no available servers for type $i$,
the customer is blocked and no further action is taken. If the customer is assigned to a server,
the server availability state ($I$ bits) is changed to indicate the unavailability to {\em any} customer type $i$.\\
(b) Each server, when its configuration changes,
i.e. upon any customer arrival (assignment) or departure (service completion),
sends a ``pull-message" ($I$ bits), containing its new availability state, to the router.\\
(c) When router receives a pull-message from a server, it updates its availability status accordingly.
(In reality, to prevent router from using ``obsolete" pull-messages, after assigning a customer
to a server, router can use some short time-out for the server, during which 
the server is considered unavailable regardless of its availability state. Thus, when 
 the time-out expires, the availability state of the server is that from the {\em latest}
pull-message received from it. If the time-out is longer than
the ``round-trip'' router-server-router message delay, then the latest pull-message from 
the server is generated upon the last customer assignment to it, or maybe later, upon departures that occurred after that.)
\\This mechanism is such that the rate of pull-messages in the system is very small, namely two pull-messages
per each arriving customer. The low rate of communication between the router and the servers is 
a very important feature of pull-based algorithms, 
because
in modern cloud based systems, the number of servers can be very large.
\\We also note that a key part of the PULL algorithm is the random uniform assignment of customers to available servers. Therefore, GRAND-F algorithm can be viewed as an extension of PULL algorithm 
to service systems with packing constraints.
\end{remark}

We consider the asymptotic regime, where the arrival rates are increased linearly 
with a scaling parameter $r\to\infty$:  $\Lambda_i = \lambda_i r$, where $\lambda_i>0$ are fixed parameters.
In addition, so do the server pool sizes $H_s$, namely,
$H_s = h_s r$, where $h_s>0, ~s\in\cs,$ are fixed parameters.

Let $\bX^r(\cdot)$ be the process associated with a system with parameter $r$, and 
let $\bX^r(\infty)$ be the (random) system state in the stationary regime.
For each $i$, denote by $Y^r_i(t) \equiv \sum_{\bk\in\ck} k_i X^r_{\bk}(t)$ the total number
of customers of type $i$. As mentioned above,  
 $Y^r_i(\infty)$ is stochastically dominated by a Poisson random variable with mean $r \rho_i$, where $\rho_i\equiv \lambda_i/\mu_i$. 
As before, without loss of generality, we assume $\sum_i \rho_i=1$.

The {\em fluid-scaled} process is $\bx^r(t)=\bX^r(t)/r$, $t \in [0, \infty)$. 
We define $\bx^r(\infty) = \bX^r(\infty)/r$.
Similarly, $y^r_i(t)=Y^r_i(t)/r$,  $x^r_{\bZero^s}(t)=X^r_{\bZero^s}(t)/r$ and 
$x^r_{(i)}(t)=X^r_{(i)}(t)/r$,  for  $t \ge 0$ and $t=\infty$.

For any $r$, $\bx^r(t)$ takes values in the compact set 
$$
\cx^{\Box} \equiv \{\bx\in \R_+^{|\ck|} ~|~  \sum_{\bk\in \ck^s} x_{\bk} \le h_s, ~\forall s\in\cs \}.
$$
For any $\bx\in\cx^{\Box}$, we denote $x_{\bZero^s}\equiv h_s - \sum_{\bk\in \ck^s} x_{\bk}, ~s\in\cs$,
and will sometimes use notation $\bar\bx \equiv \{x_{\bk}, ~\bk\in\bar\ck\}$.

The sequence of distributions of $\bx^r(\infty)$ is obviously tight,
and therefore there always exists a limit $\bx(\infty)$ in distribution, so that $\bx^r(\infty)\implies \bx(\infty)$,
along a subsequence of $r$. 
The limit (random) vector $\bx(\infty)$ satisfies the following property w.p.1.:
\beql{eq-cons-laws-555}
\sum_{\bk\in\ck} k_i x_{\bk}(\infty) \equiv y_i(\infty) \le \rho_i, ~~\forall i.
\end{equation}
The asymptotic regime and property \eqn{eq-cons-laws-555} obviously hold for any placement algorithm,
not just GRAND-F.

Consider the following subset of $\cx^{\Box}$:
$$
\cx^{\diamond} \equiv \{\cx\in \cx^{\Box} ~|~  \sum_s \sum_{\bk\in \ck^s} k_i x_{\bk} = \rho_i, ~\forall i\in\ci \} \equiv \cx^{\Box} \cap \cx.
$$
We make the following
\begin{assumption}
\label{assum-feasible}
The system parameters $\lambda_i$, $\mu_i$, $i\in\ci$, and $h_s$, $s\in\cs$, are such that the set
$\cx^{\diamond}$ in non-empty. Moreover,
there exists $\bx\in\cx^{\diamond}$ such that $x_{\bZero^s}>0$ for all $s$.
\end{assumption}

This assumption means that, when the scaling parameter $r$ is large, 
and we have $\rho_i r$
customers of each type $i$, it is possible to ``pack" all of them into the system servers 
($h_s r$ for each type $s$), so that a non-zero fraction of servers in each pool $s$ 
remains idle. Recall that, when $r$ is large, $\rho_i r$ is essentially the maximum
number of type $i$ customers the system can possibly have in steady state, because this would be
the number of customers in the infinite-server system with no blocking. Thus, the assumption
guarantees that it is  feasible, at least in principle,
 to operate a system in a way such that,
in the $r\to\infty$ limit, the steady-state blocking probability vanishes. 

Consider the following function $L^{\Box}(\bar\bx)$ defined on $\bar\bx$ such that 
$\bx\in\cx^{\Box}$ (and $x_{\bZero^s} \equiv h_s - \sum_{\bk\in\ck^s} h_{\bk}$ for all $s$):
\beql{eq-L-def-finite}
L^{\Box}(\bar\bx) = \sum_{\bk\in\bar\ck} x_{\bk} \log [x_{\bk} c_{\bk} /e],
\end{equation}
where $c_{\bk} \doteq \prod_i k_i !$, $0!=1$. 
We then have 
\beql{eq-L-partial-finite}
(\partial / \partial x_{\bk}) L^{\Box}(\bar\bx) =\log [x_{\bk} c_{\bk}], ~~\bk\in \bar\ck.
\end{equation}
For each $\bk\in\bar\ck$ the corresponding summand in the definition \eqn{eq-L-def-finite} of function 
$L^{\Box}(\bar\bx)$ is strictly convex in $x_{\bk}$; then, $L^{\Box}(\bar\bx)$ is strictly convex
on $\R_+^{|\bar\ck|}$.

Consider the problem $\min_{\bx\in \cx^\diamond} L^{\Box}(\bar\bx)$. It is the
 following convex optimization problem:
\beql{eq-opt-grand-finite}
\min_{\bar\bx\in\R_+^{|\ck|}} L^{\Box}(\bar\bx),
\end{equation}
subject to
\beql{eq-cons-laws-grand-finite}
\sum_{\bk\in\ck} k_i x_{\bk} = \rho_i, ~~\forall i,
\end{equation}
\beql{eq-cons-laws-grand555-finite}
\sum_{\bk\in\bar\ck^s} x_{\bk} = h_s, ~s\in\cs.
\end{equation}

Denote by $\bar\bx^{*,\Box}$ its unique optimal solution;
of course, the corresponding $\bx^{*,\Box} \in \cx^\diamond$.
Using \eqn{eq-L-partial-finite} and Assumption~\ref{assum-feasible} it is easy to see that $x^{*,\Box}_{\bk}>0$ for all $\bk \in \bar\ck$.
There exist a vector of Lagrange multipliers $\bnu^{*,\Box}=(\nu_i^{*,\Box}, ~ i\in \ci)$
for 
the constraints \eqn{eq-cons-laws-grand-finite} and Lagrange multipliers $\beta_s^{*}$ 
for the
constraints \eqn{eq-cons-laws-grand555-finite}, such that $\bar\bx^{*,\Box}$ solves problem
$$
\min_{\bar\bx\in \R_+^{|\bar\ck|}} L^{\Box}(\bar\bx) + \sum_i \nu_i^{*,\Box} (\rho_i - \sum_{\bk\in\ck} k_i x_{\bk})
+ \sum_s \beta_s^{*} (\sum_{\bk\in\bar\ck^s} x_{\bk} - h_s).
$$
We see that $\log [x_{\bk}^{*,\Box} c_{\bk} ] - \sum_i \nu_i^{*,\Box} k_i +\beta_s^{*} = 0$, $ \bk \in \bar\ck^s$. 
Therefore, $\bar\bx^{*,\Box}$ has the product form
\beql{eq-grand-product-finite}
x^{*,\Box}_{\bk} = \frac{1}{c_{\bk}} \exp \left[-\beta_s^{*} + \sum_i k_i \nu_i^{*,\Box}\right] = \frac{e^{-\beta_s^{*}}}{c_{\bk}} \exp \left[ \sum_i k_i \nu_i^{*,\Box}\right]  , ~~ \bk \in \bar \ck^s.
\end{equation}
This in particular implies that Lagrange multipliers $\nu_i^{*,\Box}$, $\beta_s^{*}$, are unique.
They
can have any sign (not necessarily non-negative). 

We obtain the following fact. {\em
A point $\bar\bx$, such that $\bx\in \cx^\diamond$, is the optimal solution to \eqn{eq-opt-grand-finite}-\eqn{eq-cons-laws-grand555-finite}
(that is $\bar\bx=\bar\bx^{*,\Box}$) if and only if it has a product form representation
\eqn{eq-grand-product-finite} for some Lagrange multipliers $\nu_i^{*,\Box}$, $\beta_s^{*}$.} 
Furthermore, {\em $\bx^{*,\Box}$ and $\bnu^{*,\Box}$ are equal to $\bx^{*,\ba}$ 
and $\bnu^{*,\ba}$, respectively, defined for the
infinite-server system in Section~\ref{sec-results-infinite}, with parameters $a_s=e^{-\beta_s^{*}}$.} 

Our main result for the finite-server system is the following Proposition~\ref{prop-main-finite}.
(It is stated here informally. Formal statements are given in Lemmas~\ref{lem-fsp-conv-local-on-diamond} 
and \ref{lem-fsp-conv-local}.)

\begin{prop}
\label{prop-main-finite}
Suppose Assumption~\ref{assum-feasible} holds. 
As $r\to\infty$, the limits of the fluid-scaled trajectories $\bx^r(\cdot)$ will be referred to as fluid sample paths (FSP). Point $\bx\in\cx^{\Box}$ is an invariant point, if $\bx(t)\equiv \bx$ is an FSP.
Then $\bx^{*,\Box}$ is the unique invariant point $\bx$,
such that $x_{\bZero^s}>0$ for all $s$ (and therefore there is no blocking).
Moreover, this invariant point is locally stable: $\bx(t)\to \bx^{*,\Box}$, uniformly for all FSPs with 
$\bx(0)$ sufficiently close to $\bx^{*,\Box}$.
\end{prop}

In turn, Proposition~\ref{prop-main-finite} strongly suggests that the following asymptotic optimality property holds,
which we present as

\begin{conjecture}
\label{th-grand-fluid-finite}
Suppose Assumption~\ref{assum-feasible} holds.
Consider a sequence of systems under the GRAND-F algorithm, indexed by $r$,
and let $\bx^r(\infty)$ denote the random state of the fluid-scaled process
in the stationary regime.
Then, as $r\to\infty$,
$
\bx^r(\infty) \implies \bx^{*,\Box}.
$
\end{conjecture}

If Conjecture~\ref{th-grand-fluid-finite} is correct, the GRAND-F algorithm is asymptotically optimal in the following sense. As long as Assumption~\ref{assum-feasible} holds, 
i.e. the system has enough capacity to process all offered load (under ideal packing), then
as $r\to\infty$, the steady-state blocking probability under GRAND-F vanishes.
As discussed in Remark~\ref{rem-pull}, GRAND-F can be viewed as an extension of 
PULL algorithm \cite{St2014_pull}. Therefore, Conjecture~\ref{th-grand-fluid-finite}, if correct,
can be viewed as an extension (to systems with packing constraints) of the asymptotic optimality of PULL.

\section{Proof of Theorem~\ref{th-grand-fluid-convergence}}
\label{sec-opt-conv}

For any $\bk\in\ck^s$, as $a_s\downarrow 0$,
$$
[-\log a_s]^{-1}  x_{\bk} \log [x_{\bk} c_{\bk} /(e a_s)] - x_{\bk}
= [-\log a_s]^{-1}  x_{\bk} [\log x_{\bk} + \log c_{\bk} -1] \to 0,
$$
uniformly on any compact subset of non-negative $x_{\bk}$.
We have 
$$
L^{(\ba)}(\bx)/[-\log a_1] = \sum_s [-\log a_s]/[-\log a_1]  
\sum_{\bk\in\ck^s} [-\log a_s]^{-1} x_{\bk} \log [x_{\bk} c_{\bk} /(e a_s)].
$$
Setting $a_s = \alpha^{\gamma_s}$ (which implies $[-\log a_s]/[-\log a_1]= \gamma_s/\gamma_1=\gamma_s$), we see that, as $\alpha\downarrow 0$,
 $|L^{(\ba)}(\bx)/[-\log \alpha] - \sum_s \sum_{\bk\in \ck^s} \gamma_s x_{\bk}| \to 0$, uniformly
in $\bx\in \cx$. Therefore, $\bx^{*,\ba}$ must converge to $\cx^*$.

Consider any sequence $\alpha\downarrow 0$. We will denote $b=-\log \alpha$.
We will show that from any subsequence we can choose a further subsequence, along which
we have convergence $\bx^{*,\ba} \to \bx^*$, $\bnu^{*,\ba}/b \to \bbeta^*$,
 where $\bx^* \in \cx^*$ and $\bbeta^* \in \ch^*$ .

Let a subsequence of $\alpha$ be fixed. Since  $\bx^{*,\ba} \to \cx^*$,
we can and do choose a further subsequence along which $\bx^{*,\ba} \to \bx^*$ 
for some fixed $\bx^{*} \in \cx^*$.
Let us show that 
\beql{eq-limsup-1}
\limsup_{\alpha\to 0} \sum_i k_i  \nu^{*,\ba}_i/b \le \gamma_s, ~~\forall \bk\in\ck^s,
\end{equation}
\beql{eq-liminf-0}
\liminf_{\alpha\to 0} \nu^{*,\ba}_i/b \ge 0, ~~\forall i.
\end{equation}
From \eqn{eq-grand-product} we have:
\beql{eq-grand-product777}
x^{*,\ba}_{\bk} = \frac{1}{c_{\bk}} \exp \left[b(\sum_i k_i \nu_i^{*,\ba}/b - \gamma_s)\right], ~~ \bk \in \ck^s.
\end{equation}
If \eqn{eq-limsup-1} would not hold for some $\bk\in\ck^s$, then by \eqn{eq-grand-product777}
we would have $\limsup x^{*,\ba}_{\bk} = \infty$ -- a contradiction. Thus, \eqn{eq-limsup-1} holds.
Suppose now that \eqn{eq-liminf-0} does not hold for some $i$, that is 
$\liminf \nu^{*,\ba}_i/b < 0$. Pick an $s$ and $\bk\in\ck^s$ such that $k_i \ge 1$ and $x^*_{\bk}>0$.
Such $s$ and $\bk$ must exist, because $\sum_{\bk} k_i x^*_{\bk} = \rho_i$ (recall that 
$\bx^* \in \cx^*$). Since $x^{*,\ba}_{\bk} \to x^*_{\bk} \in [0,\rho_i]$,
we see from \eqn{eq-grand-product777}
that $\lim \sum_j k_j  \nu^{*,\ba}_j/b = \gamma^s$. Therefore, 
$$
\limsup \left[\sum_{j\ne i} k_j  \nu^{*,\ba}_j/b +
(k_i - 1) \nu^{*,\ba}_i/b\right] = \gamma^s- \liminf  \nu^{*,\ba}_i/b > \gamma^s;
$$
 but, this violates \eqn{eq-limsup-1}
for configuration $\bk-\be_i$. Thus, \eqn{eq-liminf-0} holds.

By \eqn{eq-limsup-1}-\eqn{eq-liminf-0}, the sequence of $\bnu^{*,\ba}/b$ is bounded.
Then, we choose a further subsequence along which $\bnu^{*,\ba}/b$ converges to some $\bbeta^*$.
For the pair $\bx^*$ and $\bbeta^*$, condition \eqn{eq-dual-1} is automatic,
conditions \eqn{eq-dual-111}-\eqn{eq-dual-2} follow from
\eqn{eq-limsup-1}-\eqn{eq-liminf-0}, and condition \eqn{eq-dual-3} follows
from \eqn{eq-grand-product777}. Therefore, $\bbeta^* \in \ch^*$. $\Box$

\section{Fluid sample paths for the infinite-server system\\
under GRAND($\ba Z$). Proof of Theorem~\ref{th-grand-fluid} 
}
\label{sec-fsp-dynamics}

In this section, we define fluid sample paths (FSP) 
for the system controlled by GRAND($\ba Z$).
FSPs 
arise as limits of the (fluid-scaled) trajectories $(1/r) \bX^r(\cdot)$
as $r\to\infty$. Then we prove Theorem~\ref{th-grand-fluid}.
The development in this section is a generalization to the heterogeneous system of the definitions and results
 given for the homogeneous system in Section 4 of \cite{StZh2013}.
The generalization is quite straightforward. However, we provide it here for completeness
and, more importantly, 
as a preparation for the related argument used later in Section~\ref{sec-fsp-dynamics-finite}
 for the finite-server system.

Let $\cm$ denote the set of pairs $(\bk,i)$ such that $\bk\in\ck$ and $\bk-\be_i\in\bar\ck$.
Each pair $(\bk,i)$ is associated with the ``edge'' $(\bk-\be_i,\bk)$ connecting configurations
$\bk-\be_i$ and $\bk$; often we refer to this edge as $(\bk,i)$.
By ``arrival along the edge $(\bk,i)$'', we will mean placement of a type $i$ customer
into a server configuration $\bk-\be_i$ to form configuration $\bk$. Similarly, ``departure along the edge $(\bk,i)$'' is a departure of a type-$i$ customer
from a server in configuration $\bk$, which changes its configuration to $\bk-\be_i$.

Without loss of generality, assume that
the Markov process $X^r(\cdot)$ for each $r$ is driven by the
common set of primitive processes, defined as follows.

For each $(\bk,i)\in \cm$, consider an independent
unit-rate Poisson process $\{\Pi_{\bk i}(t), ~t\ge 0\}$, which drives 
departures along edge $(\bk,i)$. Namely, 
let $D^r_{\bk i}(t)$ denote the total 
number of  departures along the edge $(\bk,i)$ in $[0,t]$; then
\beql{eq-driving-dep}
D^r_{\bk i}(t) = \Pi_{\bk i} \left(\int_{\bZero}^t X_{\bk}^r(s) k_i \mu_i ds\right).
\end{equation}
The functional strong law of large numbers (FSLLN) holds:
\beql{eq-flln-poisson-dep}
\frac{1}{r}\Pi_{\bk i}(rt) \to t, ~~~u.o.c., ~~w.p.1.
\end{equation}
For each $i\in \ci$, consider an independent
unit-rate Poisson process $\Pi_{i}(t), ~t\ge 0$, which drives 
exogenous arrivals of type $i$. Namely, 
let $A^r_{i}(t)$ denote the total 
number of type-$i$ arrivals in $[0,t]$, then
\beql{eq-driving-arr}
A^r_{i}(t) = \Pi_{i}(\lambda_i r t).
\end{equation}
Analogously to \eqn{eq-flln-poisson-dep},
\beql{eq-flln-poisson-arr}
\frac{1}{r}\Pi_{i}(rt) \to t, ~~~u.o.c., ~~w.p.1.
\end{equation}
The random placement of new arrivals is constructed as follows.
For each $i\in \ci$, consider an i.i.d. sequence
$\xi_i(1), \xi_i(2), \ldots$ of random variables, uniformly distributed in $[0,1]$.
Denote by $\ck_i \doteq \{\bk\in \bar\ck ~|~ \bk+\be_i \in \bar\ck\}$ the subset of those 
configurations (including zero configurations) which can fit an additional type-$i$ 
customer.
The configurations $\bk\in \ck_i$ are indexed by $1,2,\ldots,|\ck_i|$ 
(in arbitrary fixed order). When the $m$-th (in time) customer 
of type $i$ arrives in the system, it is assigned as follows.
If $X_{(i)}^r=0$, the customer is assigned to an empty server 
of an arbitrarily fixed type $s$, such that $\be_i \in \ck^s$.
Suppose $X_{(i)}^r \ge 1$. Then, the customer
is assigned 
to a server in configuration $\bk'$ indexed by $1$ if 
$$
\xi_i(m) \in [0,X^r_{\bk'}/ X^r_{(i)}],
$$
it is assigned to
a server in configuration $\bk''$ indexed by $2$ if 
$$
\xi_i(m) \in (X^r_{\bk''}/ X^r_{(i)},(X^r_{\bk'}+X^r_{\bk''})/ X^r_{(i)}],
$$
and so on.  Denote
$$
g^r_i(\sigma,\zeta) \doteq \sum_{m=1}^{\lfloor r\sigma \rfloor} I\{\xi_i(m) \le \zeta\},
$$
where $\sigma\ge 0$, $0\le \zeta \le 1$. 
Obviously, from the strong law of large numbers
 and the monotonicity of $g^r_i(\sigma,\zeta)$ on both arguments, we 
have the FSLLN
\beql{eq-flln-random}
g^r_i(\sigma,\zeta) \to \sigma\zeta, ~~~\mbox{u.o.c.} ~~~\mbox{w.p.1}
\end{equation}

It is easy (and standard) to see that, for any $r$, w.p.1,
 the realization of the process $\{\bX^r(t), ~t\ge 0\}$
is uniquely determined by the
initial state $\bX^r(0)$
and the realizations of the driving processes
$\Pi_{\bk i}(\cdot)$, $\Pi_{i}(\cdot)$ and $(\xi_i(1), \xi_i(2), \ldots)$.

If we denote by $A^r_{\bk i}(t)$ the total number of arrivals 
allocated along edge $(\bk,i)$ in $[0,t]$, we obviously have
$\sum_{\bk\in \ck_i} A^r_{\bk i}(t)= A^r_{i}(t), ~t\ge 0$, for each $i$.

In addition to
$$
x^r_{\bk}(t) = \frac{1}{r} X^r_{\bk}(t),
$$
we introduce other  fluid-scaled 
quantities:
$$
d^r_{\bk i}(t) = \frac{1}{r} D^r_{\bk i}(t),~~~
a^r_{\bk i}(t) = \frac{1}{r} A^r_{\bk i}(t).
$$

A set of locally Lipschitz continuous functions
$[\{x_{\bk}(\cdot),~\bk\in \ck\}, \{d_{\bk i}(\cdot),~(\bk,i)\in \cm\},\{a_{\bk i}(\cdot),~(\bk,i)\in \cm\}]$
on the time interval $[0,\infty)$ we call a {\em fluid sample path} (FSP), if there exist
realizations of the primitive driving processes,
 satisfying conditions \eqn{eq-flln-poisson-dep},\eqn{eq-flln-poisson-arr}
and \eqn{eq-flln-random}
and a fixed subsequence of $r$, along which
\begin{eqnarray}
& [\{x_{\bk}^r(\cdot),~\bk\in \ck\}, \{d_{\bk i}^r(\cdot),~(\bk,i)\in \cm\},\{a_{\bk i}^r(\cdot),~(\bk,i)\in \cm\}]
\to \nonumber \\
& [\{x_{\bk}(\cdot),~\bk\in \ck\}, \{d_{\bk i}(\cdot),~(\bk,i)\in \cm\},\{a_{\bk i}(\cdot),~(\bk,i)\in \cm\}],
~~u.o.c. \label{eq-fsp-def-closed}
\end{eqnarray}

For any FSP, all points $t>0$ are regular (see definition in Section~\ref{subsec-notation}),
except a subset of zero Lebesgue measure. 

\begin{lem}
\label{lem-conv-to-fsp-closed}
Consider a sequence of fluid-scaled processes $\{\bx^r(t),~t\ge 0\}$
with fixed initial states $\bx^r(0)$ such that $\bx^r(0)\to \bx(0)$.
Then w.p.1, for any subsequence of $r$ there exists a further
subsequence of $r$, along which the convergence \eqn{eq-fsp-def-closed} holds,
with the limit being an FSP.
\end{lem}

{\em Proof} is same as that of Lemma 5 in \cite{StZh2013}. 
$\Box$

For an FSP, at a regular time point $t$, we denote
$v_{\bk i}(t)=(d/dt)a_{\bk i}(t)$ and $w_{\bk i}(t)=(d/dt)d_{\bk i}(t)$.
In other words, $v_{\bk i}(t)$ and $w_{\bk i}(t)$ are the rates of type-$i$ ``fluid''
arrival and departure along edge $(\bk,i)$, respectively. Also
denote: $y_i(t)=\sum_{\bk} k_i x_{\bk}(t)$, $z(t)=\sum_i y_i(t)$, $x_{\bZero^s}(t)=a_s z(t)$, 
and $x_{(i)}(t) = \sum_{\bk \in \bar\ck: \bk + \be_i \in \bar\ck} x_{\bk}(t)$.

\begin{lem}
\label{lem-fsp-properties-closed-basic}
(i) An FSP satisfies the following properties
at any regular point $t$:
\beql{eq-y-ode}
(d/dt) y_i(t) = \lambda_i - \mu_i y_i(t), ~~\forall i\in \ci,
\end{equation}
\beql{eq-grand-dep-rate}
w_{\bk i}(t)=k_i \mu_i x_{\bk}(t), ~~\forall (\bk,i)\in \cm,
\end{equation}
\beql{eq-grand-arr-rate}
x_{(i)}(t)>0 ~~\mbox{implies}~~
v_{\bk i}(t)=\frac{x_{\bk-\be_i}(t)}{x_{(i)}(t)} \lambda_i, ~~\forall (\bk,i)\in \cm,
\end{equation}
\beql{eq-conserv-rate}
\sum_{\bk:(\bk,i)\in \cm} v_{\bk i}(t) = \lambda_i, ~~\forall i\in \ci,
\end{equation}
\beql{eq-main-diff}
(d/dt) x_{\bk}(t) = \left[\sum_{i:\bk-\be_i\in\bar\ck} v_{\bk i}(t) - \sum_{i:\bk+\be_i\in\bar\ck} v_{\bk+\be_i,i}(t) \right]
                     - \left[\sum_{i:\bk-\be_i\in\bar\ck} w_{\bk i}(t) - \sum_{i:\bk+\be_i\in\bar\ck} w_{\bk+\be_i,i}(t) \right], ~~\forall \bk\in\ck.
\end{equation}
Clearly, \eqn{eq-y-ode} implies
\beql{eq-y-ode2}
y_i(t) = \rho_i + (y_i(0)-\rho_i) e^{-\mu_i t}, ~~t\ge 0, ~~\forall i\in \ci.
\end{equation}
(ii) Moreover, an FSP with $\bx(0)\in \cx$ satisfies the following stronger
conditions:
\beql{eq-y-ode-star}
y_i(t) \equiv \rho_i, ~~\forall i\in \ci,
\end{equation}
\beql{eq-a}
z(t) \equiv 1, ~~ x_{\bZero^s}(t)\equiv a_s, ~~ x_{(i)}(t) \ge \sum_{s:~\be_i\in\ck^s} a_s, ~\forall i\in \ci;
\end{equation}
at any regular point $t$,
\beql{eq-grand-arr-rate-star}
v_{\bk i}(t)=\frac{x_{\bk-\be_i}(t)}{x_{(i)}(t)} \lambda_i, ~~\forall (\bk,i)\in \cm,
\end{equation}
\beql{eq-conserv-rate-star}
\sum_{\bk:(\bk,i)\in \cm} w_{\bk i}(t) = 
\lambda_i, ~~\forall i\in \ci.
\end{equation}
\end{lem}

{\em Proof.}  (i) Given the convergence
\eqn{eq-fsp-def-closed}, which defines an FSP, all the stated properties
except \eqn{eq-grand-arr-rate} are
nothing  but the limit versions of the flow conservations laws.
Property \eqn{eq-grand-arr-rate} follows from the construction 
of the random assignment, the continuity of $\bx(t)$, and \eqn{eq-flln-random}.
We omit further details. \\
(ii) If $\bx(0)\in \cx$, which implies $y_i(0)=\rho_i$ for each $i$,
property \eqn{eq-y-ode-star} (and then \eqn{eq-a} as well) 
follows from \eqn{eq-y-ode2}. Then, \eqn{eq-grand-arr-rate}
strengthens to \eqn{eq-grand-arr-rate-star}, and 
\eqn{eq-conserv-rate-star} is verified directly using \eqn{eq-grand-dep-rate}.
$\Box$

\begin{lem}
\label{lem-fsp-conv}
For any FSP with $\bx(0)\in \cx$,
\beql{eq-fsp-conv}
\bx(t) \to \bx^{*,\ba},
\end{equation}
and the convergence is uniform across all such FSPs.
\end{lem}

{\em Proof.} 
Given that $x_{\bZero^s}(t) \equiv a_s$ and $\sum_{\bk} x_{\bk}(t) \le 1$, we have 
$x_{(i)}(t) \le 1+\sum_s a_s$, 
hence $v_{\bk i}(t) \ge x_{\bk}(t) \lambda_i / (1+\sum_s a_s)$.
From here, we obtain the following fact: for any $\bk$ and
any $\delta>0$ there exists $\delta_1>0$ such that for all $t\ge \delta$,
 $x_{\bk}(t) \ge \delta_1$. The proof is by 
contradiction. Consider a $\bk$, say $\bk\in\bar\ck^s$, that is a minimal counterexample; 
necessarily, $\bk \ne \bZero^s$.
Pick any $\delta>0$ and then the corresponding $\delta_1>0$ such that 
the statement holds for any $\bk'\in \bar\ck^s$, $\bk'<\bk$. 
(Here $\bk'<\bk$ means that $\bk'_i\le \bk_i, ~\forall i,$ and $\bk'\ne \bk$.)
We observe from \eqn{eq-main-diff} that
 for any regular $t\ge \delta$, $(d/dt) x_{\bk}(t) > \delta_2 >0$ as long as
$x_{\bk}(t) \le \delta_3$, for some positive constants $\delta_2, \delta_3$.
Since this holds for an arbitrarily small $\delta>0$ (with 
$\delta_1, \delta_2, \delta_3$ depending on it), we see that the statement
is true for $\bk$.

In particular, we see that $x_{\bk}(t)>0$ for all $t>0$ and all $\bk$.
Note also that all $t>0$ are regular points (because all $w_{\bk i}$
and $v_{\bk i}$ are bounded continuous in $\bx$).

To prove the lemma, it will suffice to show that:\\
 (a) if $\bx(t) \ne  \bx^{*,\ba}$ 
and $x_{\bk}(t)>0$ 
for all $\bk \in \ck$, then $(d/dt) L^{(\ba)}(\bx(t))<0$; and, moreover,\\
(b) the derivative is bounded away from zero as long as $\| \bx(t) -  \bx^{*,\ba}\|$ 
is bounded away from zero.\\
Let us denote by $\Xi(\bx)$ the derivative $(d/dt) L^{(\ba)}(\bx(t))$ at a given point
$\bx(t)=\bx$; in the rest of the proof we study the function $\Xi(\bx)$ on $\cx$, 
and therefore drop the time index $t$. Suppose all components $x_{\bk}>0$.
From \eqn{eq-grand-dep-rate}, \eqn{eq-conserv-rate},
\eqn{eq-grand-arr-rate-star}, and
\eqn{eq-conserv-rate-star}, we have:
\beql{eq-w-expr}
w_{\bk i} = k_i \mu_i x_{\bk} = k_i \mu_i x_{\bk}\sum_{\bk': (\bk',i) \in \cm} \frac{x_{\bk'-\be_i}}{x_{(i)}},
\end{equation}
\beql{eq-v-expr}
v_{\bk' i} = \frac{x_{\bk'-\be_i}}{x_{(i)}} \lambda_i = \frac{x_{\bk'-\be_i}}{x_{(i)}} \sum_{\bk: (\bk,i) \in \cm} k_i \mu_i x_{\bk}.
\end{equation}
Expressions \eqn{eq-w-expr} and \eqn{eq-v-expr} can be interpreted as follows.
For any ordered pair of edges $(\bk,i)$ and
$(\bk',i)$, we can assume that the part $k_i \mu_i x_{\bk} x_{\bk'-\be_i}/x_{(i)}$
of the total departure rate $k_i \mu_i x_{\bk}$ along $(\bk,i)$ is ``allocated back''
as a part of the arrival rate along $(\bk',i)$. 
Using \eqn{eq-L-partial},
the contribution of these ``coupled'' departure/arrival
rates for the ordered pair of edges $(\bk,i)$ and
$(\bk',i)$ into the derivative $\Xi(\bx)$ is 
$$
\xi_{\bk,\bk',i} = \left[\log (k'_i x_{\bk-\be_i} x_{\bk'}) - \log (k_i x_{\bk} x_{\bk'-\be_i}) \right]
\frac{k_i \mu_i x_{\bk} x_{\bk'-\be_i}}{x_{(i)}}.
$$
This expression is valid even when either $\bk-\be_i=\bZero^s$ or $\bk'-\be_i=\bZero^s$
for some $s$.
This is because 
$x_{\bZero^s}(t) = a_s$ when $\bx\in \cx$, and therefore by convention \eqn{eq-formal-deriv-zero},
formula \eqn{eq-L-partial} is valid for all $\bk\in \bar\ck$. We have:
$$
\xi_{\bk,\bk',i} + \xi_{\bk',\bk,i} 
= (\mu_i/x_{(i)}) [\log (k'_i x_{\bk-\be_i} x_{\bk'}) - \log (k_i x_{\bk} x_{\bk'-\be_i}) ] 
[k_i x_{\bk} x_{\bk'-\be_i} - k'_i x_{\bk-\be_i} x_{\bk'}] \le 0,
$$
and the inequality is strict unless $k'_i x_{\bk-\be_i} x_{\bk'}=k_i x_{\bk} x_{\bk'-\be_i}$.
We obtain
\beql{eq-L-deriv}
\Xi(\bx) = \sum_i \sum_{\bk,\bk'} [\xi_{\bk,\bk',i} + \xi_{\bk',\bk,i}].
\end{equation}
Therefore, 
$\Xi(\bx)<0$ unless $\bx$ has a product form representation \eqn{eq-grand-product},
which in turn is equivalent to $\bx=\bx^{*,a}$.

So far the function $\Xi(\bx)$ in \eqn{eq-L-deriv} was defined for $\bx\in\cx$
with all $x_{\bk}>0$. Let us adopt a convention that $\Xi(\bx)=-\infty$
for $\bx\in\cx$ with at least one $x_{\bk}=0$.  Then,
it is easy to verify that $\Xi(\bx)$ is continuous on the entire set $\cx$.

It remains to show that
for any $\delta_2 > 0$ there exists $\delta_3 > 0$ such that conditions
$\bx\in \cx$ and
$L^{(\ba)}(\bx) - L^{(\ba)}(\bx^{*,\ba}) \ge \delta_2$ imply
$\Xi(\bx) \le -\delta_3$. This is indeed true, because otherwise
there would exist $\bx\in \cx$, $\bx\ne \bx^{*,\ba}$, such that 
 $\Xi(\bx) = 0$, which is, again, equivalent to $\bx = \bx^{*,\ba}$.
$\Box$

From Lemma~\ref{lem-fsp-conv} we easily obtain Theorem~\ref{th-grand-fluid};
see the proof of Theorem 3 in Section 4 of \cite{StZh2013}.

As in \cite{StZh2013}, we also have the following generalization of Lemma~\ref{lem-fsp-conv},
showing FSP uniform convergence  for
arbitrary initial states, not necessarily $\bx(0)\in \cx$. 

\begin{lem}
\label{lem-fsp-conv-gen}
For any compact $A \in \R_+^{|\ck|}$, the convergence
\beql{eq-fsp-conv-gen}
\bx(t) \to \bx^{*,\ba}
\end{equation}
holds uniformly in all FSPs with $\bx(0)\in A$.
\end{lem}

{\em Proof} repeats that of Lemma 8 in \cite{StZh2013} almost verbatim.
The only adjustments are: \\
1) Starting any fixed time $\tau>0$, we have $0< a_1 \le x_{\bZero^s}(t), ~\forall s,$
and $x_{(i)}(t) \le a_2 < \infty, ~\forall i,$ for some constants $a_1,a_2$, uniformly 
on all FSPs with $\bx(0)\in A$;\\
2) $L^{(\ba)}$ replaces $L^{(a)}$;\\
3) $f(\bk) = (\partial / \partial x_{\bk}) L^{(\ba)}(x) =\log [x_{\bk} c_{\bk} / a_s ], ~\bk\in\ck^s.$
$\Box$

\section{GRAND-F: Local stability of FSPs}
\label{sec-fsp-dynamics-finite}

The construction of the Markov process $X^r(\cdot)$ under GRAND-F is the
same as in Section~\ref{sec-fsp-dynamics} for GRAND($\ba Z$), except now, when 
$X_{(i)}^r=0$, an arriving type $i$ customer is blocked.
Consequently, we no longer have
the identity
$\sum_{\bk\in \ck_i} A^r_{\bk i}(t)= A^r_{i}(t), ~t\ge 0$, for each $i$. Instead,
$$
A^r_{i}(t) - \sum_{\bk\in \ck_i} A^r_{\bk i}(t), ~t\ge 0,
$$
is non-negative non-decreasing function, giving the number of blocked type $i$ customers by time $t$.

The definition of an FSP and Lemma~\ref{lem-conv-to-fsp-closed} hold as is.
All points $t>0$ are regular, except for a subset of zero Lebesgue measure. 
The analog of Lemma~\ref{lem-fsp-properties-closed-basic} is the following 

\begin{lem}
\label{lem-fsp-properties-closed-basic-fin}
(i) An FSP satisfies the following properties
at any regular point $t$:
\beql{eq-conserv-rate-fin}
\sum_{\bk:(\bk,i)\in \cm} v_{\bk i}(t) \le \lambda_i, ~~\forall i\in \ci,
\end{equation}
\beql{eq-y-ode-fin}    
(d/dt) y_i(t) = \sum_{\bk:(\bk,i)\in \cm} v_{\bk i}(t) - \mu_i y_i(t), ~~\forall i\in \ci,
\end{equation}
\beql{eq-grand-dep-rate-fin}
w_{\bk i}(t)=k_i \mu_i x_{\bk}(t), ~~\forall (\bk,i)\in \cm,
\end{equation}
\beql{eq-grand-arr-rate-fin}
x_{(i)}(t)>0 ~~~\mbox{implies}~~~ \sum_{\bk:(\bk,i)\in \cm} v_{\bk i}(t) = \lambda_i, ~\forall i\in \ci,
~~\mbox{and}~~
v_{\bk i}(t)=\frac{x_{\bk-\be_i}(t)}{x_{(i)}(t)} \lambda_i, ~\forall (\bk,i)\in \cm,
\end{equation}
\beql{eq-main-diff-fin}
(d/dt) x_{\bk}(t) = \left[\sum_{i:\bk-\be_i\in\bar\ck} v_{\bk i}(t) - \sum_{i:\bk+\be_i\in\bar\ck} v_{\bk+\be_i,i}(t) \right]
                     - \left[\sum_{i:\bk-\be_i\in\bar\ck} w_{\bk i}(t) - \sum_{i:\bk+\be_i\in\bar\ck} w_{\bk+\be_i,i}(t) \right], ~~\forall \bk\in\bar\ck.
\end{equation}
(ii) Moreover, an FSP with $\bx(0)\in \cx^{\diamond}$, $x_{(i)}(0)>0,~\forall i,$ satisfies the following stronger
conditions for all sufficiently small $t>0$:
\beql{eq-y-ode-star-fin}
y_i(t) \equiv \rho_i, ~~\forall i\in \ci,
\end{equation}
\beql{eq-a-fin}
z(t) \equiv 1, ~~ x_{\bZero^s}(t)\equiv a_s,~\forall s, ~~ x_{(i)}(t) \ge \min_s a_s, ~\forall i\in \ci;
\end{equation}
if $t$ is regular,
\beql{eq-grand-arr-rate-star-fin}
v_{\bk i}(t)=\frac{x_{\bk-\be_i}(t)}{x_{(i)}(t)} \lambda_i, ~~\forall (\bk,i)\in \cm,
\end{equation}
\beql{eq-conserv-rate-star-fin}
\sum_{\bk:(\bk,i)\in \cm} w_{\bk i}(t) = 
\lambda_i, ~~\forall i\in \ci.
\end{equation}
\end{lem}

{\em Proof.} 
(i) Given the convergence
\eqn{eq-fsp-def-closed} defining an FSP, all the stated properties
except \eqn{eq-grand-arr-rate-fin}, are
nothing  but the limit versions of the flow conservations laws.
Property \eqn{eq-grand-arr-rate-fin} follows from the construction 
of the random assignment, the continuity of $\bx(t)$, and \eqn{eq-flln-random}.
We omit further details. \\
(ii) If $\bx(0)\in \cx^{\diamond}$, which implies $y_i(0)=\rho_i$ for each $i$,
property \eqn{eq-y-ode-star-fin} (and then \eqn{eq-a-fin} as well) 
follows from \eqn{eq-y-ode-fin}  and \eqn{eq-grand-arr-rate-fin}. 
Then, \eqn{eq-conserv-rate-star-fin} is verified directly using \eqn{eq-grand-dep-rate-fin}.
Finally, \eqn{eq-grand-arr-rate-star-fin} follows from
\eqn{eq-grand-arr-rate-fin}.
$\Box$

\begin{lem}
\label{lem-fsp-conv-local-on-diamond}
There exists $\epsilon>0$, such that, uniformly on FSPs with initial states $\bx(0) \in \cx^\diamond \cap \{\|\bx - \bx^{*,\Box}\| \le \epsilon\}$, 
\beql{eq-fsp-convlocal-on-diamond}
\bx(t) \to \bx^{*,\Box}, ~~t\to\infty.
\end{equation}
FSP $\bx(t)\equiv \bx^{*,\Box}$ is the unique invariant FSP, satisfying conditions 
$x_{\bZero^s}(0) > 0, ~\forall s$.
\end{lem}

{\em Proof.} We can assume (without loss of generality) that $\epsilon$ is small enough so that
$x_{\bk}(0)>0, ~\forall \bk\in\bar\ck$. In particular, at $t=0$, the condition 
$x_{\bZero^s}(t)>0, ~\forall s$, holds. Obviously, until the first time $\tau>0$ when this condition is violated
($\tau=\infty$ if it is never violated), we have $y_i(t) = \rho_i, ~\forall i$. It is also easy to see that all time ponts $0< t < \tau$ are regular and such that $x_{\bk}(t)>0, ~\forall \bk\in\bar\ck$. 
Denote by $\Xi(\bar\bx)$ the derivative $(d/dt) L^{\Box}(\bar\bx(t))$ at a given point
$\bx(t)=\bx$. Then, expressions \eqn{eq-w-expr} and \eqn{eq-v-expr} for $w_{\bk i}$ and 
$v_{\bk' i}$ hold for our system, and can be interpreted the same way.
(Recall, however, that now the components $x_{\bZero^s}$ are {\em not} constant,
and therefore their derivatives do depend on the rates $w_{\bZero^s+\be_i, i}$ and 
$v_{\bZero^s+\be_i, i}$.)
Then the expression for $\Xi(\bar\bx)$ has exactly same form as expression \eqn{eq-L-deriv} for $\Xi(\bx)$
in Section~\ref{sec-fsp-dynamics}:
\beql{eq-L-deriv-fin}
\Xi(\bar\bx) = \sum_i \sum_{(\bk,i),(\bk',i)} = 
(\mu_i/x_{(i)}) [\log (k'_i x_{\bk-\be_i} x_{\bk'}) - \log (k_i x_{\bk} x_{\bk'-\be_i}) ] 
[k_i x_{\bk} x_{\bk'-\be_i} - k'_i x_{\bk-\be_i} x_{\bk'}] \le 0.
\end{equation}
The inequality in \eqn{eq-L-deriv-fin}
is strict unless $k'_i x_{\bk-\be_i} x_{\bk'}=k_i x_{\bk} x_{\bk'-\be_i}$
for all pairs of edges $(\bk,i)$ and $(\bk',i)$.
Therefore, 
$\Xi(\bar\bx)<0$ unless $\bar\bx$ has a product form representation \eqn{eq-grand-product-finite},
which in turn is equivalent to $\bx=\bx^{*,\Box}$. 

Function $\Xi(\bar\bx)$ is continuous in a neighborhood of  $\bx^{*,\Box}$
(and in fact at any point such that $x_{\bk} > 0, ~\forall \bk\in\bar\ck$).
Choose $\epsilon_1>0$ small enough so that $x_{\bk}>0, ~ \bk\in\bar\ck,$ for all 
$\bx \in \cx^\diamond \cap \{\|\bx - \bx^{*,\Box}\| \le \epsilon_1\}$.
Then choose $\delta>0$ such that condition $L^{\Box}(\bar\bx) - L^{\Box}(\bar\bx^{*,\Box}) \le \delta$
(along with $\bx \in \cx^\diamond$) implies $\|\bx - \bx^{*,\Box}\| < \epsilon_1$.
Finally, choose $\epsilon>0$ small enough so that 
the maximum of $L^{\Box}(\bar\bx) - L^{\Box}(\bar\bx^{*,\Box})$
over the set $\cx^\diamond \cap \{\|\bx - \bx^{*,\Box}\| \le \epsilon\}$ is less than $\delta$.
We see that a trajectory with $\bx(0) \in \cx^\diamond \cap \{\|\bx - \bx^{*,\Box}\| \le \epsilon\}$
cannot escape from the set $\cx^\diamond \cap \{\|\bx - \bx^{*,\Box}\| \le \epsilon_1\}$,
and therefore $x_{\bk}(t)>0, ~\bk\in\bar\ck,$ for all $t\ge 0$. Then, the convergence 
\eqn{eq-fsp-convlocal-on-diamond} holds, and it is uniform on $\bx(0) \in \cx^\diamond \cap \{\|\bx - \bx^{*,\Box}\| \le \epsilon\}$, because, for any $0<\delta_1 < \delta$, $\Xi(\bar\bx)$  is negative and bounded away from zero for all $\bx\in\cx^\diamond \cap \{\delta_1 \le L^{\Box}(\bar\bx) - L^{\Box}(\bar\bx^{*,\Box}) \le \delta\}$.

It is a corollary from the above argument, that there cannot be an invariant FSP
$\bx(t)\equiv \bx(0)$ with
$x_{\bZero^s}(0) > 0, ~\forall s$, unless $\bx(0) = \bx^{*,\Box}$.
(Indeed, $\bx(0) \in \cx^\diamond$ necessarily, because if $y_i(0)\ne \rho_i$ then $y_i(t)$ cannot be constant.
Then $\bx(0) = \bx^{*,\Box}$, because
otherwise $L^{\Box}(\bar\bx(t))$ cannot be constant.)
This proves the second statement of the lemma. $\Box$

\begin{lem}
\label{lem-fsp-conv-local}
There exists $\epsilon>0$, such that, uniformly on FSPs with initial states $\bx(0) \in \cx^\Box \cap \{\|\bx - \bx^{*,\Box}\| \le \epsilon\}$, 
\beql{eq-fsp-convlocal-on-diamond222}
\bx(t) \to \bx^{*,\Box}, ~~t\to\infty.
\end{equation}
\end{lem}

{\em Proof} is a slightly generalized version of that of Lemma~\ref{lem-fsp-conv-local-on-diamond}.
That proof considers FSPs that stay within $\cx^\diamond$, uses the continuity of 
$\Xi(\bar\bx)$, and the fact that for $\bx\in \cx^\diamond$ in a small
neighborhood of $\bx^{*,\Box}$, $\Xi(\bar\bx)<0$  unless $\bx=\bx^{*,\Box}$.
But, $\Xi(\bar\bx)$ is continuous in a neighborhood of  $\bx^{*,\Box}$
(or any point such that $x_{\bk} > 0, ~\forall \bk\in\bar\ck$), not necessarily restricted to $\cx^\diamond$.
In addition, we know that as long as $x_{\bZero^s}(t)>0, ~\forall s$,
each $y_i(t)$ satisfies ODE 
$
(d/dt) [y_i(t) - \rho_i] = -\mu_i [y_i(t) - \rho_i],
$
and therefore 
\beql{eq-conv-to-diamond}
y_i(t) - \rho_i = (y_i(0) - \rho_i)e^{-\mu_i t}.
\end{equation}

Using these observations, the adjustment of the proof of Lemma~\ref{lem-fsp-conv-local-on-diamond}
is as follows. We choose small  $\epsilon_1>0$, then $\delta>0$, then $\epsilon>0$, exactly as in that proof.
Then, using the continuity of $\Xi(\bar\bx)$, along with \eqn{eq-conv-to-diamond},
we can choose a sufficiently small $\epsilon_2>0$, so that
a trajectory with $\bx(0) \in \{|y_i-\rho_i| \le \epsilon_2, ~\forall i\} \cap \{\|\bx - \bx^{*,\Box}\| \le \epsilon\}$
cannot escape from the set $\{|y_i-\rho_i| \le \epsilon_2, ~\forall i\} \cap \{\|\bx - \bx^{*,\Box}\| \le \epsilon_1\}$.
Then, the convergence 
\eqn{eq-fsp-convlocal-on-diamond222} holds, and it is uniform on 
$\bx(0) \in \{|y_i-\rho_i| \le \epsilon_2, ~\forall i\} \cap \{\|\bx - \bx^{*,\Box}\| \le \epsilon\}$, 
because, for any $0<\delta_1 < \delta$, 
there exists a small $\epsilon'_2>0$, such that 
$\Xi(\bar\bx)$  is negative and bounded away from zero for all 
$\bx\in \{|y_i-\rho_i| \le \epsilon'_2, ~\forall i\} \cap \{\delta_1 \le L^{\Box}(\bar\bx) - L^{\Box}(\bar\bx^{*,\Box}) \le \delta\}$. (Note that the time for FSPs starting in 
$\{|y_i-\rho_i| \le \epsilon_2, ~\forall i\} \cap \{\|\bx - \bx^{*,\Box}\| \le \epsilon\}$
to reach set $\{|y_i-\rho_i| \le \epsilon'_2, ~\forall i\} \cap \{\|\bx - \bx^{*,\Box}\| \le \epsilon\}$
is uniformly bounded due to \eqn{eq-conv-to-diamond}.)
$\Box$

\subsection{Comments on Conjecture~\ref{th-grand-fluid-finite}, local stability, and fixed point argument}
\label{sec-fixed-point}

Lemmas~\ref{lem-fsp-conv-local-on-diamond} and \ref{lem-fsp-conv-local} formally state properties 
described informally in Proposition~\ref{prop-main-finite}. The sequence of steady-states $\bx^r(\infty)$
 is obviously tight. 
It is easy to see that its any subsequential limit in distribution, $\bx(\infty)$, is such that
$y_i(\infty) \le \rho_i, ~\forall i,$ w.p.1. This is because, by comparison with the infinite-server system,
 $Y_i^r(\infty)$ is stochastically dominated
by a Poisson random variable with mean $\rho_i r$. Furthermore, again by comparison with the infinite-server system, any FSP with
$$
\bx(0) \in \cx^{\Box,\le} \equiv \{\bx\in \cx^{\Box} ~|~  \sum_s \sum_{\bk\in \ck^s} k_i x_{\bk} \le \rho_i, ~\forall i\in\ci \}
$$
stays in $\cx^{\Box,\le}$ at all times $t$. Given these facts, if we would have the 
(analogous to Lemma~\ref{lem-fsp-conv}) uniform convergence
property 
\beql{eq-conv-desire}
\bx(t) \to \bx^{*,\Box}, ~~~\forall \bx(0)\in \cx^{\Box,\le},
\end{equation}
this would prove Conjecture~\ref{th-grand-fluid-finite} (by the same argument as in the proof 
of Theorem~\ref{th-grand-fluid}). Unfortunately, the uniform convergence \eqn{eq-conv-desire}
does {\em not} hold for a general finite-server system. It is very easy to construct a counterexample
(e.g., for a system with one server type with the configuration set shown of Fig. 1(b) in \cite{StZh2013})
such that there exists an invariant FSP $\bx(t)\equiv \bx^*$, ``sitting" at a suboptimal point
$\bx^* \ne \bx^{*,\Box}$, such that $y_i^* < \rho_i, ~\forall i$, and therefore such that there is
non-zero fraction of customers of each type being blocked. (In fact, we believe that a stronger property
holds for such a counterexample: the sequence of processes $\bx^r(\cdot)$ converges {\em in distribution}
to the invariant FSP $\bx(t)\equiv \bx^*$.) This, of course, does not imply that Conjecture~\ref{th-grand-fluid-finite} is wrong -- it just shows that there is no hope of proving Conjecture~\ref{th-grand-fluid-finite}
based on fluid scale considerations alone.

Lemmas~\ref{lem-fsp-conv-local-on-diamond} and \ref{lem-fsp-conv-local} show FSP local stability at the optimal point  $\bx^{*,\Box}$, and the fact that $\bx^{*,\Box}$ is the only invariant point at which there is no blocking.
This strongly suggests that Conjecture~\ref{th-grand-fluid-finite} is correct, even though, 
as discussed above, it is insufficient for its proof. Still, we note that the local stability is a substantially
stronger property than a typical ``fixed point" argument which is used to ``guess" asymptotic properties
like our Conjecture~\ref{th-grand-fluid-finite}. In our case a ``fixed point" argument would go as follows:
as $r\to\infty$, assume that steady-state distributions of server states are asymptotically independent;
further assume that a subsequential limit of the marginal distribution of a server state is such that
the server is empty with non-zero probability; under these assumptions, find the set of (limiting) marginal 
distributions (for each server type), which would remain invariant (``fixed") over time; 
in our case, this argument leads to finding that the only such possible set of marginal distributions
is such that the system must be ``sitting" at the point $\bx^{*,\Box}$, equal to the one defined in this paper. 
Note that, in essence, the above argument
is nothing else but the statement that $\bx^{*,\Box}$ is the unique invariant point (at which there is no blocking)
for FSPs, while local stability properties in Lemmas~\ref{lem-fsp-conv-local-on-diamond} and \ref{lem-fsp-conv-local} are much stronger.

\section{Discussion}
\label{sec-further-work}

Proving Conjecture~\ref{th-grand-fluid-finite} for the finite-server system under GRAND-F is a very interesting and challenging
subject of future work. As discussed in Section~\ref{sec-fixed-point}, fluid-scale analysis alone 
cannot  be sufficient for such a proof, because there may exists sub-optimal points, which are invariant for the FSPs.

The local stability results for the finite-server system with blocking
(Proposition~\ref{prop-main-finite}, Lemmas~\ref{lem-fsp-conv-local-on-diamond} 
and \ref{lem-fsp-conv-local}) hold for other variants of the finite-server system as well.
Indeed, these results and their proofs 
only concern with the system behavior in the vicinity of equilibrium point, where there are always
available servers for any customer type. Suppose now that we have a system in which customers are queued instead of blocking
when there are no available servers for them
(or a system where both blocking and queueing are possible). Then the local stability results still apply
for this system, {\em as long as 
the assignment rule coincides with GRAND-F when there are servers available to arrivals.}
Further, this suggests that Conjecture~\ref{th-grand-fluid-finite} is also valid 
for such other variants of the finite-server system, under appropiate versions of GRAND-F.
In fact, recall that GRAND-F, as defined in this paper,
itself can be
viewed as an extension of PULL algorithm \cite{St2014_pull} to systems with packing constraints. 
PULL algorithm has been defined and proved to be asymptotically optimal 
for very general systems with queueing and/or blocking (but without packing constraints).

The results of this paper further highlight the universality of GRAND algorithm.
For example, Best Fit type algorithms are applicable only to the special case of vector packing constraints, where the underlying notion of a customer ``fitting best into the remaining space'' at a server makes sense. When packing constraints are more general, Best Fit is not applicable, while GRAND is.
Furthermore, inherently, Best Fit requires precise information about the current state of each server -- this can be a substantial disadvantage in practical large-scale systems. GRAND, on the other hand, only needs to know whether a given customer fits into a given server or not; this allows a very efficient practical implementation (as discussed in detail in  Remark~\ref{rem-pull}). It is possible that 
versions of Best Fit may perform better than GRAND for systems with vector packing constraints. Paper \cite{GZS2013} provides some evidence of that. 
(Although, the algorithm studied in \cite{GZS2013} is not a ``pure'' Best Fit,
but a Best Fit {\em with randomization}, a mixture, in a sense, of Best Fit and GRAND.) 
Studying versions of Best Fit is an interesting subject; it is outside the scope of this paper, which is focused on general packing constraints. First Fit is another approach to packing;
algorithms of this type use fixed preordering of servers and place each customer into the first one where it can fit.
Such algorithms are easily implementable and apply to general packing constraints. Note that
GRAND can be viewed as a First Fit with random uniform reordering of servers before each customer placement.
If the order of servers has to chosen and fixed a priori, as ``pure'' First Fit requires,
the question arises on how to do it when the servers
are heterogeneous, as in our model. Exploring variants of First Fit may be another subject of future research.

\end{document}